\documentclass[11pt]{article}

\usepackage{latexsym}

\usepackage{amsmath}
\usepackage{amssymb}

\usepackage{pst-all}

\input epsf

\def\d#1{{#1\kern-0.4em\char"16\kern-0.1em}}
\def\D#1{{\raise0.2ex\hbox{-}\kern-0.4em  #1}}
\def\Dj{\text{\raise0.3ex\hbox{-}\kern-0.4em  D}}

\def\be{\begin{equation}}
\def\ee{\end{equation}}
\def\ba{\begin{eqnarray*}}
\def\ea{\end{eqnarray*}}
\def\baa{\begin{eqnarray}}
\def\eaa{\end{eqnarray}}

\def\e{\varepsilon}

\def\a{\alpha}
\def\F{{\varPhi}}
\def\la{\lambda}

\def\bm{\boldmath}

\title{\Large\bf On root-ratio multipoint methods for finding multiple zeros of univariate functions}

\author{{\bf Miodrag S. Petkovi\'c$^{\,1,}$, Ljiljana D. Petkovi\'c$^{\,2,}$}\footnote{Corresponding author}
   \\[2mm]
   $^1$\small\it  Faculty of Electronic Engineering,
University of Ni\v s,\\[-1mm]
\small\it A. Medvedeva 14, 18000 Ni\v s, Serbia\\
$^2$\small\it  Faculty of Mechanical Engineering,
University of Ni\v s, \\[-1mm] \small \it A. Medvedeva 14, 18000 Ni\v s, Serbia}

\date{}

\begin{document}

\maketitle

\begin{abstract}{\small Several root-ratio multipoint methods for finding multiple zeros of univariate functions were recently presented. The characteristic of these methods is that they deal with $m$-th root of ratio of two functions (hence the name root-ratio methods), where $m$ is the multiplicity of the sought zero, known in advance. Some of these methods were presented without any derivation and motivation, it could be said, out of the blue. In this paper we present an easy and entirely natural way for constructing root-ratio multipoint iterative methods starting from  multipoint methods for finding simple zeros. In this way, a vast number of root-ratio multipoint methods for  multiple zeros, existing as well new ones, can be constructed. For demonstration, we derive four root-ratio methods for multiple zeros. Besides, we study computational cost of the considered methods and give a comparative analysis that involves CPU time needed for the extraction of the $m$-th root. This analysis shows  that root-ratio methods are pretty inefficient from the computational point of view and, thus, not suitable in practice.  A general discussion on a practical need for multipoint methods of very high order is also considered.
}
\\[3mm]
{\small{\it AMS Mathematical Subject Classification (2010):}
65H05.
\\ [1mm]
}
\end{abstract}

\renewcommand{\thefootnote}{}%
\footnote{{\it E-mail address:}  ljiljana.petkovic@masfak.ni.ac.rs (L. D. Petkovi\'c)
 }

\date{}

\maketitle

\section{Introduction}

In the last ten years a lot of papers were published in the topic of iterative methods of optimal order for finding multiple zeros of univariate functions. Among them, several methods were constructed using {\it root-parameter approach} dealing with parameters computed by extracting $m$-th root of some values. In this paper this class of methods will be called {\it root-ratio methods}. We  concentrate on the two issuses: (i) computational cost of these methods and (ii)  demonstration of an easy procedure for natural construction of root-ratio methods starting from methods for finding simple zeros.

\smallskip

The paper is organized as follows. In Section 2 we present three two-point root-ratio methods of optimal order four published in \cite{zhou-AMC}, \cite{lee-AMC} and \cite{liu}.  In Section 3 we expose three tree-point root-ratio methods of optimal order eight proposed in \cite{behl}, \cite{zafar-1} and \cite{zafar-3}.  An easy and fully natural procedure for constructing root-ratio methods of optimal order, demonstrated by three examples, is presented in Section 4. The computational  efficiency of root-ratio methods is given in Section 5 applying a comparative analysis that involves CPU time needed for the extraction of the $m$-th root. This analysis has shown that root-ratio methods are pretty expensive and not-competitive with multipoint methods which do not require the extraction of the root even if the latter ones have lower order of convergence. See Remark 8 at the end of the paper.

\section{Two-point fourth order methods}

Let $\a$ be a zero of the known multiplicity  $m\ge 1$ of a differentiable function $f,$ and $k=0,1,\ldots$ be the iteration index. We will assume that the chosen initial approximation $x_0$ is sufficiently close to the zero $\a$ so that the considered iterative processes are convergent.
We start with two-point root-ratio methods for finding multiple zero of functions of one variable. Zhou et al. proposed in \cite{zhou-AMC} the following family of iterative two-step methods for finding multiple zeros:
\be
\left\{\begin{array}{l}
y_k=x_k-m\dfrac{f(x_k)}{f'(x_k)},\\[9pt]
x_{k+1}=y_k-mG(u_k) \dfrac{f(x_k)}{f'(x_k)},\quad u_k=\Bigl(\dfrac{f(y_k)}{f(x_k)}\Bigr)^{1/m}\label{1}
\end{array}
\right.\quad (k=0,1,\ldots).
\ee
Under the conditions $G(0)=0,\ G'(0)=1,\ G''(0)=4,\ G'''(0)<+\infty,$ the order of convergence of the methods (\ref{1}) is four. The iterative formula was derived correctly but in a complicated way.

\smallskip

Lee et al. constructed in  \cite{lee-AMC} the following family of fourth order methods:
\be
\left\{\begin{array}{l}
y_k=x_k-m\cdot \dfrac{f(x_k)}{f'(x_k)+\la f(x_k)},\\[13pt]
x_{k+1}=y_k-m W(u_k)\cdot \dfrac{f(x_k)}{f'(x_n)+2\la f(x_k)}, \quad u_k=\Bigl(\dfrac{f(y_k)}{f(x_k)}\Bigr)^{1/m},\label{2}
\end{array}
\right.
\ee
where
$$
W(u)=\frac{u(1+(c+2)u+ru^2)}{1+cu}
$$
and  $\la,\ c,\ r$ are arbitrary parameters.

\bigskip

{\bf Remark 1.}  The parameter $\la$ in (\ref{2}) does not improve convergence characteristics of the methods (\ref{2}). Indeed, if $m\ge 2$, then both values $f'(x_k)$ and $f(x_k)$ are very small in magnitude if the approximation $x_k$ is sufficiently close to the zero $\a$ of $f$ and, consequently,  the parameter $\la$ does not play any important role so that one can take $\la=0$ in (\ref{2}) without loss of generality. The derivatives of $W$ are given by
$$
W'(u)=\frac{1+2(c+2)u+(2c+c^2+3r)u^2+2cru^3}{(1+cu)^2},\quad W''(u)=\frac{4+2ru(3+3cu+c^2u^2)}{(1+cu)^3}.
$$
Setting $u=0$ one obtains $W(0)=0,\ \ W'(0)=1,\ \ W''(0)=4.$ Therefore, the method (\ref{2}) is essentially a special case of the more general method (\ref{1}).

\bigskip
Using root-ratio approach in a slightly different way,  Liu et al. derived in \cite{liu} the family of  two-point  iterative methods
\be
\left\{\begin{array}{l}
y_k=x_k-m\dfrac{f(x_k)}{f'(x_k)},\\[9pt]
x_{k+1}=y_k-mG(t_k) \dfrac{f(x_k)}{f'(x_k)},\quad t_k=\Bigl(\dfrac{f'(y_k)}{f'(x_k)}\Bigr)^{\frac{1}{m-1}}\label{3}
\end{array}
\right. \quad (k=0,1,\ldots),
\ee
which converges with the order of convergence four under the conditions
$$
G(0)=0,\quad G'(0)=1,\quad G''(0)=\frac{4m}{m-1}.
$$
The family (\ref{3})  cannot be applied for finding simple zeros ($m=1$), which is an obvious disadvantage. This drawback limits its applications in composite algorithms where a zero-finding method is  part of the algorithm in which the case $m=1$ is possible; actually,  this case often appears in practical problems (for example, solving engineering problems).

\smallskip

All three two-point methods (\ref{1}), (\ref{2}) and (\ref{3}) are {\it optimal} in the sense of Kung-Traub conjecture \cite{KT}: {\it The highest order of convergence of an $n$-point method requiring $n+1$ function evaluations is $2^n$.}

\section{Three-point methods of eighth order}

In this section we present some recent three-point optimal iterative methods for finding multiple zeros of the known multiplicity.

\smallskip

Behl et al.  proposed in \cite{behl} (without derivation and motivation) the family of  three-point methods of optimal order eight:
\be
\left\{\begin{array}{l}
y_k=x_k-m\dfrac{f(x_k)}{f'(x_k)},\\[9pt]
z_k=y_k-u_k S(h_k)  \dfrac{f(x_k)}{f'(x_k)},\\[9pt]
x_{k+1}=z_k-u_kv_k R(h_k,v_k) \dfrac{f(x_k)}{f'(x_k)} \label{4}
\end{array}
\right. \quad (k=0,1,\ldots),
\ee
where
$$
u_k=\Bigl(\frac{f(y_k)}{f(x_k)}\Bigr)^{1/m},\quad h_k=\frac{u_k}{a_1+a_2u_k}\ (a_1\ne 0),\quad v_k=\Bigl(\frac{f(z_k)}{f(y_k)}\Bigr)^{1/m}.
$$
Under some specific conditions for the function $S$ (at the point 0) and $R$ (at the point (0,0)), the order of convergence of the family of iterative methods (\ref{4}) is eight.

\smallskip

Zafar et al. presented  in \cite{zafar-1} (without derivation and motivation) the following iterative three-point method
\be
\left\{\begin{array}{l}
y_k=x_k-m\dfrac{f(x_k)}{f'(x_k)},\\[9pt]
z_k=y_k-mu_k H(u_k)  \dfrac{f(x_k)}{f'(x_k)},\\[9pt]
x_{k+1}=z_k-u_kv_k(B_1+B_2u_k) P(v_k)G(w_k) \dfrac{f(x_k)}{f'(x_k)} \label{5}
\end{array}
\right. \quad (k=0,1,\ldots),
\ee
where $B_1,B_2\in \mbox{\bm{$R$}}$ are free parameters,  $H,\ P,\ G$ are the weight functions and
\be
u_k=\Bigl(\frac{f(y_k)}{f(x_k)}\Bigr)^{1/m},\quad
v_k=\Bigl(\frac{f(z_k)}{f(y_k)}\Bigr)^{1/m},\quad \label{6}
w_k=\Bigl(\frac{f(z_k)}{f(x_k)}\Bigr)^{1/m}.
\ee

Zafar et al. proposed  in \cite{zafar-3} (without derivation and motivation)  the family of three-point methods
\be
\left\{\begin{array}{l}
y_k=x_k-m\dfrac{f(x_k)}{f'(x_k)},\\[9pt]
z_k=y_k-mu_k H(u_k)  \dfrac{f(x_k)}{f'(x_k)},\\[9pt]
x_{k+1}=z_k-u_kP(u_k)G(v_k)L(w_k) \dfrac{f(x_k)}{f'(x_k)} \label{7}
\end{array}
\right.\quad (k=0,1,\ldots),
\ee
where $u,\ v,\ w$ are defined by (\ref{6}). The family (\ref{7}) is obvious   generalization of the family (\ref{5}). Indeed, comparing (\ref{5}) and (\ref{7}), we note that $P(u)=B_1+B_2u,\ G(v)=vP(v).$ Clearly,  the authors could work immediately at start with the product of weight functions $P(u)G(v)L(w)$ since the generalization is quite obvious.

\bigskip

{\bf Remark 2.} Neither  the idea nor the motivation for the construction of the methods (\ref{4}), (\ref{5}) and (\ref{7}) were presented; the  iterative formulas appear {\it out of the blue}, without preliminary explanation/introduction of the basic idea and derivation procedures. Such non-preamble approach is not in the spirit of the methodology of scientific work and educational principles; hence it is of a little significance for readers.
\bigskip

{\bf Remark 3.} The need for arbitrary parameters in any iterative formula is discussed in Remark 7.  This relates to the  parameters $\la,c,r$ in (\ref{2}), and especially to $a_1$ and $a_2$
in (\ref{4}). Since any advantage of using $a_1$ and $a_2$ was not proved in \cite{behl}, it is logical to take $a_1=1$ and $a_2=0$ without loss of generality. At  first sight, this is a special case but, in fact, it  is a  natural  choice which simplifies the iterative formula (\ref{4}). Consequently, it follows $h_k=u_k$ and the family (\ref{4}) reduces to a simpler iterative method (without $h_k$)
\be
\left\{\begin{array}{l}
y_k=x_k-m\dfrac{f(x_k)}{f'(x_k)},\\[9pt]
z_k=y_k-u_k S(u_k)  \dfrac{f(x_k)}{f'(x_k)},\\[9pt]
x_{k+1}=z_k-u_kv_k R(u_k,v_k) \dfrac{f(x_k)}{f'(x_k)} \label{8}
\end{array}
\right. \quad (k=0,1,\ldots).
\ee

\section{Can we construct them easier?}

In this section we show that some of the above-presented  methods can be derived using an easy way, well known in the literature for almost 150 years. The idea was known to German mathematician E. Schr\"oder, see his paper \cite{sreder} and the English translation \cite{stjuart}, and goes as follows in the case of Newton's  method
$$
x_{k+1}=x_k-\frac{f(x_k)}{f'(x_k)}\quad (k=0,1,\ldots)
$$ applied to the function $f$ having a zero $\a$ of multiplicity $m\ge 1.$

\bigskip
4.1 {\it Schr\"oder's $f^{1/m}$-approach}
\bigskip

Let $F(x)=f(x)^{1/m}.$ Then $\a$ is a simple zero of the function $F.$  Applying Newton's method to $F$, we obtain
$$
x_{k+1}=x_k-\frac{F(x_k)}{F'(x_k)}=x_k-\frac{f(x)^{1/m}}{\frac{1}{m}f'(x)f(x)^{1/m-1}}
$$
and hence
$$
x_{k+1}=x_k-m\frac{f(x_k)}{f'(x_k)} \quad (k=0,1,\ldots),
$$
which is well-known Schr\"oder's method for finding multiple zero of the known multiplicity $m.$ This useful idea was applied in many papers during the last 70 years.

\smallskip

Schr\"oder's approach can be applied to multipoint methods for finding a simple zero in order to construct corresponding multipoint mehods for multiple zeros. As above, let $f$ be a function having a zero of the known multiplicity $m\ge 1.$ First, we set $F(x)=f(x)^{1/m}$ (Schr\"oder's $f^{1/m}$-approach) and derive the following relations:
\baa
\frac{F(x)}{F'(x)}&=&\frac{f(x)^{1/m}}{\frac{1}{m}f'(x)f(x)^{1/m-1}}=m\frac{f(x)}{f'(x)},
\nonumber\\ \frac{F(y)}{F'(x)}&=& \frac{f(y)^{1/m}}{\frac{1}{m}f'(x)f(x)^{1/m-1}}=m \frac{f(x)}{f'(x)}\Bigl(\frac{f(y)}{f(x)}\Bigr)^{1/m}=m u\frac{f(x)}{f'(x)},\quad u=\Bigl(\frac{f(y)}{f(x)}\Bigr)^{1/m},\nonumber\\
\frac{F(z)}{F'(x)}&=& \frac{f(z)^{1/m}\cdot \frac{f(y)^{1/m}}{f(y)^{1/m}}}{\frac{1}{m}f'(x)f(x)^{1/m-1}}=m \frac{f(x)}{f'(x)}\Bigl(\frac{f(y)}{f(x)}\Bigr)^{1/m}\Bigl(\frac{f(z)}{f(y)}\Bigr)^{1/m}=m u v \frac{f(x)}{f'(x)},\quad v=\Bigl(\frac{f(z)}{f(y)}\Bigr)^{1/m},\nonumber\\
\frac{F(y)}{F(x)}&=&\frac{f(y)^{1/m}}{f(x)^{1/m}}=u,\label{9}\\ \frac{F(z)}{F(y)}&=&\frac{f(z)^{1/m}}{f(y)^{1/m}}=v,\nonumber\\ \frac{F(z)}{F(x)}&=&\frac{f(z)^{1/m}}{f(x)^{1/m}}=w=uv.\nonumber
\eaa

\medskip

4.2 {\it Eighth-order family {\rm (\ref{4})} derived by Schr\"oder's approach}

\bigskip
For demonstration, we will apply the presented Schr\"oder's approach  to the
 following three-point family of iterative methods for finding a simple zero of a function  $f,$  proposed in \cite{AMC-2011}:
\be
\left\{\begin{array}{l}
y_k=x_k-\dfrac{f(x_k)}{f'(x_k)},\\[9pt]
z_k=y_k-p(u_k) \dfrac{f(y_k)}{f'(x_k)},\\[9pt]
x_{k+1}=z_k-q(u_k,v_k) \dfrac{f(z_k)}{f'(x_k)} \label{10}
\end{array}
\right.\quad (k=0,1,\ldots),
\ee
where
$$
u_k=\frac{f(y_k)}{f(x_k)},\quad v_k=\frac{f(z_k)}{f(y_k)}
$$
and $p(u)$ and $q(u,v) $ are functions of one and two variables, respectively. The following theorem has been proved  in \cite{AMC-2011} (see, also, \cite{elzevir}).

\bigskip

{\bf Theorem 1.} {\it Let $a,\ b$ and $c$ be arbitrary constants. If $p$ and $q$ are arbitrary differentiable functions with Taylor's series of the form
\ba
p(u)&=&1+2u+\frac{a}{2}u^2+\frac{b}{6}u^3+\cdots,\\
q(u,v)&=&1+2u+v+\frac{2+a}{2}u^2+4uv+\frac{c}{2}v^2+\frac{6a+b-24}{6}u^3+\cdots,\\
\ea
then the family of three-point methods {\rm (\ref{10})} is of order eight. It is assumed that higher-order terms are represented by the dots, and they can take arbitrary values.}

\medskip

Assume that $f$ has the zero $\a$ of multiplicity $m\ge 1,$ known in advance. We modify (\ref{10}) by using the above-mentioned  Schr\"oder's strategy with the function
$
F(x)=f(x)^{1/m}
$
for which $\a$ is a simple zero. First, replacing $f$ with $F$, rewrite (\ref{10}) in the form
\be
\left\{\begin{array}{l}
y_k=x_k-\dfrac{F(x_k)}{F'(x_k)},\\[9pt]
z_k=y_k-P(u_k) \dfrac{F(y_k)}{F'(x_k)},\\[9pt]
x_{k+1}=z_k-Q(u_k,v_k) \dfrac{F(z_k)}{F'(x_k)} \label{11}
\end{array}
\right.\quad (k=0,1,\ldots),
\ee
where
$$
u_k=\frac{F(y_k)}{F(x_k)},\quad v_k=\frac{F(z_k)}{F(y_k)}.
$$

Using the relations (\ref{9}), the family (\ref{11}) can be written in the following form:
\be
\left\{\begin{array}{l}
y_k=x_k-m \dfrac{f(x_k)}{f'(x_k)},\\[9pt]
z_k=y_k-m u_kP(u_k) \dfrac{f(x_k)}{f'(x_k)},\\[9pt]
x_{k+1}=\F(x_k):=z_k-m u_kv_kQ(u_k,v_k) \dfrac{f(x_k)}{f'(x_k)} \label{12}
\end{array}
\right. \quad (k=0,1,\ldots),
\ee
where
\be
u_k=\Bigl(\frac{f(y_k)}{f(x_k)}\Bigr)^{1/m},\quad
v_k=\Bigl(\frac{f(z_k)}{f(y_k)}\Bigr)^{1/m}.\quad \label{13}
\ee
The weight functions $P(u)$ and $Q(u,v)$ should be determined in such a way that the iterative methods defined by (\ref{12}) have the maximal order of convergence using only the calculation of $f(x),\;f'(x),\;f(y)$ and $f(z).$ In the ideal case,  the order would reach {\it eight}, and the family (\ref{12}) would be {\it optimal} according to Kung-Traub's hypothesis \cite{KT}.

\smallskip

We proceed using a standard technique based on Taylor's series by employing symbolic computation in computer algebra system {\it Mathematica}. Since this technique was seen many times in existing papers, we present only an outline of convergence analysis.

\smallskip

We omit the iteration index $k$ and define the errors
$$
\e=x-\a,\ \ \e_y=y-\a,\ \ \e_z=z-\a,\ \ \hat \e=\hat x-\a,
$$
where $\hat x$ is a new approximation $x_{k+1}.$
Introduce
$$
 C_r=\frac{m!}{(m+r)!} \frac{f^{(m+r)}(\a)}{f^{(m)}(\a)}\ \ (r=1,2,\ldots).
    $$
We will use the following development of the function $f$ about the zero $\a$ of multiplicity $m$
    $$
    f(x)=\frac{f^{(m)}(\a)}{m!}\e^m\Bigl(1+C_1\e+C_2\e^2+C_3\e^3+C_4\e^4+C_5\e^5+
    C_6\e^6+C_7\e^7+C_8\e^8+O(\e^9)\Bigr),
    $$
    and  a program in  {\it Mathematica}. As usual, in finding the weight functions $P$ and $Q,$ we represent these functions by their Taylor's series at the neighborhood of $u=0$ (for $P$) and $(u,v)=(0,0)$ (for $Q$):
    \ba
    &&\hspace{-0.6cm}P(u)=P(0)+P'(0)u+\dfrac{P''(0)}{2}u^2+\cdots,\\
    &&\hspace{-0.6cm}Q(u,v)=Q(0,0)+Q_{u0}(0,0)u+Q_{0v}(0,0)v+
    \frac{1}{2!}\Bigl(Q_{uu}(0,0)u^2+2Q_{uv}(0,0)uv+Q_{vv}(0,0)v^2\Bigr)+\cdots
    \ea
  Subscript indices denote partial derivatives.
    In the last developments, as well as in the program,
    the following notation is used:
    \ba
    &&fam=f^{(m)}(\a),\ fx=f(x),\ fy=f(y),\ fz=f(z),\ f1x=f'(x),\\
    && e=\e,\ ey=\e_y,\ ez=\e_z,\ e1=\hat{\e},\\
     && Pr=\frac{P^{(r)}(u)}{du}\Big|_{(u=0)}\ (r=0,1,2),\\
     &&Q00=Q(0,0),\ Qu0=\frac{\partial\,Q}{\partial\,u}\Bigl|_{(u,v)=(0,0)},\ Q0v=\frac{\partial\,Q}{\partial\,v}\Bigl|_{(u,v)=(0,0)},\\
     &&
     Quu=\frac{\partial^{\,2}Q}{\partial\,u^2}\Bigl|_{(u,v)=(0,0)},\
     Quv=\frac{\partial^{\,2}Q}{\partial\,u\,\partial\,v}\Bigl|_{(u,v)=(0,0)},\
     Qvv=\frac{\partial^{\,2}Q}{\partial\,v^2}\Bigl|_{(u,v)=(0,0)}.\
    \ea

    \medskip

    The coefficients of Taylor's developments of the weight functions $P$ and $Q$ are determined using an interactive approach by combining the program realized in  {\it Mathematica} (two parts) and the annihilation of  coefficients standing at $\e$ of lower power.

    \bigskip

    PART I ({\it Mathematica})

    \medskip


 {\tt   {\small
\noindent
fxx=1+C1*e+C2*e\^{\hspace{0.2mm}}2+C3*e\^{\hspace{0.2mm}}3+ C4*e\^{\hspace{0.2mm}}4+C5*e\^{\hspace{0.2mm}}5+C6*e\^{\hspace{0.2mm}}6+
    C7*e\^{\hspace{0.2mm}}7+C8*e\^{\hspace{0.2mm}}8;\\
    fx=fam/m!*e\^{\hspace{0.2mm}}m*fxx;
    fx1=D[fx,e]; newt=Series[fx/fx1,\{e,0,8\}];\\
     ey=e-m*newt;
 fyy=1+C1*ey+C2*ey\^{\hspace{0.2mm}}2+3*ey\^{\hspace{0.2mm}}3+C4*ey\^{\hspace{0.2mm}}4;\\
   yx=fyy*Series[1/fxx,\{e,0,8\}];\\
 yxm=Series[yx\^{\hspace{0.2mm}}{1/m},\{e,0,8\}]; \\
u=yxm*ey/e;\\
P=P0+P1*u+P2/2*u\^{\hspace{0.2mm}}2+P3/6*u\^{\hspace{0.2mm}}3+P4/24*u\^{\hspace{0.2mm}}4;\\ ez=Series[ey-m*u*P*newt//FullSimplify,\{e,0,8\}]
 }
 }

\medskip
This programm gives the following OUTCOME:

 \ba
\e_z&=&\sum_{r=2}^8 S_r\e^r\\
&=&\frac{(C_1 - C_1 P_0) \e^2}{m}\\
 &&+\frac{(-2 C_2 m (-1 + P_0) +
    C_1^2 (-1 + m (-1 + P_0) + 3 P_0 - P_1)) \e^3}{m^2}\\
    && +
 \frac{1}{2 m^3}\Bigl(-6 C_3 m^2 (-1 + P_0) +
    2 C_1 C_2 m (-4 + 3 m (-1 + P_0) + 11 P_0 - 4 P_1)\\
    && +
    C_1^3 (2 - 13 P_0 + 10 P_1 + m (4 - 2 m (-1 + P_0) - 11 P_0 + 4 P_1) -
       P_2)\Bigr) \e^4+\sum_{r=5}^8 U_r+O\bigl(e^9\bigr).
       \ea

       To annihilate coefficients by $\e^2$ and $\e^3$, from the condition $S_2=0$ and $S_3=0$ we find
       \be
       P_0=1,\ \ P_1=2,\ \ P_2\ \mbox{\rm arbitrary}, \label{14}
       \ee
   yielding
       \be
    \e_z=\frac{-2 mC_1 C_2  + C_1^3 (9 + m - P_2)}{2m^3} \e^4+O\bigl(\e^5\bigr).\label{14a}
    \ee

       The part II of the program uses previously found entries and serves for finding additional conditions which provide optimal order eight.

   \bigskip

    PART II - CONTINUATION ({\it Mathematica})

    \medskip
    {\tt
   { \small
 \noindent fzz=1+C1*ez+C2*ez\^{\hspace{0.2mm}}2; fyy=1+C1*ey+C2*ey\^{\hspace{0.2mm}}2+ C3*ey\^{\hspace{0.2mm}}3+C4*ey\^{\hspace{0.2mm}}4;\\
zy=fzz*Series[1/fyy\{e,0,8\}]; zym=Series[zy\^{\hspace{0.2mm}}(1/m),\{e,0,8\}];\\
v=zym*ez/ey;\\
Q=Q00+Qu0*u+Q0v*v+Quu/2*u\^{\hspace{0.2mm}}2+Qvv/2*v\^{\hspace{0.2mm}}2+Quv*u*v;\\
e1=Series[ez-m*u*v*Q*newt,\{e,0,8\}]//FullSimplify
}
 }

 \medskip

 The error $\hat \e=\hat x-\a\ (=e1)$ is given in the form
 $$
 \e_1=\sum_{r=4}^8 T_r\e^r+O\bigl(\e^9\bigr)
 $$
 From the conditions $T_4=0,\ T_5=0,\ T_6=0,\ T_7=0,$ we find  the coefficients
 \be
 Q_{00}=1,\ \, Q_{u0}=2,\ \, Q_{0v}=1,\ \, Q_{uv}=4, \ \, P_2 \ \mbox{\rm is arbitrary},\ \, P_3=24-6P_2,\ \, Q_{uu}=P_2+2.
  \label{15}
\ee

  In addition, we obtain
  \ba
T_8&=&
   -\frac{1}{48m^7}\Bigl\{C_1 \Bigr[C_1^2
   (m-P_2+9)-2m C_2\Bigr]\\
   &&\times
   \Bigl[C_1^4 \left(-14 m^2+3 Q_{tt}
   (m-P_2+9)^2
 +6 mP_2-204 m+150
   P_2-P_4-1054\right)\\
   &&-12m C_1^2
  C_2 \bigl(Q_{tt} (m-P_2+9)-4
   m+P_2-34\bigr)-24m^2C_1C_3+12m^2
   C_2^2  (Q_{tt}-2)\Bigr]\Bigr\},
   \ea
         so that we can write
      $$
      \hat \e=T_8\e^8+O(\e^9),
      $$
      that is,
      $$
     x_{k+1}-\a= \F(x_{k})-\a=O(\e_k^8),\ \ \e_k=x_k-\a.
      $$ In this way we have proved the following assertion.

      \bigskip

      {\bf Theorem 2.} {\it If the initial approximation $x_0$ is sufficiently close to the zero $\a$ of $f$ and the conditions {\rm (\ref{14})} and {(\rm \ref{15})} are valid, then the order of the three-point family {\rm (\ref{12})} is eight.}

      \medskip
    Note that $T_8$ tends to the asymptotic error constant of the family (\ref{12}) when $\e\to 0.$ According to Theorem 2, it follows that the three-point family (\ref{12}) is {\it optimal}.

\smallskip

    The simplest forms of the weight functions $P$ and $Q$ are their truncated Taylor series
    $$
    P(u)=1+2u+\beta\,u^2+(4-2\beta)\,u^3,\quad Q(u,v)=1+2u+v+4uv +(\beta+1)u^2,\quad \beta\ \mbox{\rm is arbitrary}.
    $$
Finding a list  of particular weight functions $P$ and $Q$ is a routine  work and it is left to the interested reader as exercise.

\bigskip

4.3 {\it Eighth-order family {\rm (\ref{7})}  derived by Schr\"oder's approach}

\bigskip

The presented Schr\"oder's $f^{1/m}$-approach can also be applied for the derivation of the method (\ref{7}). Let us start from the following
three-point method for simple zeros
\be
\left\{\begin{array}{l}
y_k=x_k-\dfrac{f(x_k)}{f'(x_k)},\\[9pt]
z_k=y_k-u_k H(u_k)  \dfrac{f(x_k)}{f'(x_k)}\\[9pt]
x_{k+1}=z_k-u_kP(u_k)G(v_k)L(w_k) \dfrac{f(x_k)}{f'(x_k)}
\end{array}
\right.  \quad (k=0,1,\ldots),\label{15a}
\ee
where
$$
u_k=\frac{f(y_k)}{f(x_k)},\quad
v_k=\frac{f(z_k)}{f(y_k)},\quad
w_k=\frac{f(z_k)}{f(x_k)}\, .
$$
 This method has the order 8 under the conditions
\baa
&& H(0)=1,\ H'(0)=2,\ P'(0)=2P(0),\ P''(0)=P(0)(2+H''(0)),\nonumber\\
&& L'(0)=2L(0),\ P'''(0)=P(0)(H'''(0)+6H''(0)-24),\label{cond}\\
&& G(0)=0,\ G'(0)=\dfrac{1}{L(0)P(0)},\ G''(0)=\dfrac{2}{L(0)P(0)}.\nonumber
\eaa
The denotation of weight functions in (\ref{15a}) and (\ref{cond}) is adjusted to the denotation used in  \cite{zafar-3}.

\smallskip

Applying Schr\"oder's $f^{1/m}$-approach and (\ref{9}) we obtain the family of three-point methods for finding multiple zeros
\be
\left\{\begin{array}{l}
y_k=x_k-m\dfrac{f(x_k)}{f'(x_k)},\\[9pt]
z_k=y_k-mu_k H(u_k)  \dfrac{f(x_k)}{f'(x_k)},\\[9pt]
x_{k+1}=z_k-m u_kP(u_k)G^*(v_k)L(w_k) \dfrac{f(x_k)}{f'(x_k)}
\end{array}
\right. \quad (k=0,1,\ldots),\label{15b}
\ee
which is equivalent to (\ref{7}) setting $mG^*(v)=G(v)$ (compare (\ref{7}) and (\ref{15b})).

\bigskip
{\bf Remark 4.} Among other methods, the first author of this paper presented the family of methods (\ref{15a}) in his lecture under the title {\it Multipoint methods  for solving nonlinear equations} at International conference {\it Computational Methods in Applied Mathematics} (Berlin, 2012). However, due to the similarity to the iterative formula (4.91) in the book \cite[p\;151]{elzevir},  the first author did not publish (\ref{15a}).  Other possible sources of (\ref{15a}) are not known to the authors.

\bigskip

4.4 {\it A modification of the family {\rm(\ref{7})}}

\bigskip

Since
$$
w=\Bigl(\frac{f(z)}{f(x)}\Bigr)^{1/m}=
\Bigl(\frac{f(y)}{f(x)}\Bigr)^{1/m}\Bigl(\frac{f(z)}{f(y)}\Bigr)^{1/m}=u\cdot v,
$$
the weight  function $L(w)$ can be omitted in (\ref{7}). Then we can construct in an easy way the following family of three-point methods for finding multiple zeros involving two parametric functions
\be
\left\{\begin{array}{l}
y_k=x_k-m\dfrac{f(x_k)}{f'(x_k)},\\[9pt]
z_k=y_k-mu_k H(u_k)  \dfrac{f(x_k)}{f'(x_k)},\\[9pt]
x_{k+1}=z_k-m u_k v_k(1+2u_kv_k)P(u_k)G(v_k) \dfrac{f(x_k)}{f'(x_k)}. \label{7b}
\end{array}
\right.
\ee
Making suitable changes in the above program in {\it Mathematica} (see \S\; 4.2), we prove that the method (\ref{7b}) is of order eight under the following conditions:
\ba
&&H(0)=1,\ \ H'(0)=2,\ \ H''(0)=m+9,\\
&&P(0)=1,\ \ P'(0)=2,\ \ P''(0)=m+11,\ \ P'''(0)=30+6m+H'''(0),\\
&&G(0)=1,\ \ H'(0)=1.
\ea

    \medskip

4.5 {\it Fourth order methods for multiple zeros}

\bigskip
Chun constructed in \cite{chun} the following two-point family of iterative methods for finding simple zeros
   \be
\left\{\begin{array}{l}
y_k=x_k- \dfrac{f(x_k)}{f'(x_k)},\\[9pt]
x_{k+1}=y_k-  \dfrac{f(y_k)}{h(u_k)f'(x_k)},\quad u_k=\dfrac{f(y_k)}{f(x_k)}, \label{16}
\end{array}
\right.
\ee
where $h(u)$ is  the weight function. The method (\ref{16}) is of fourth order under the condition $h(0)=1,\ h'(0)=-1, \ |h''(0)|<+\infty.$ To generate suitable two-point methods, sometimes it is necessary to develop the function $h$ into Taylor's or geometric series.

\smallskip

A slightly more direct approach without altering the weight function, based on Chun's idea, was given in \cite{petko-AADM} in the form
   \be
\left\{\begin{array}{l}
y_k=x_k- \dfrac{f(x_k)}{f'(x_k)},\\[9pt]
x_{k+1}=y_k- p(u_k) \dfrac{f(y_k)}{f'(x_k)},\quad u_k=\dfrac{f(y_k)}{f(x_k)}, \label{17}
\end{array}
\right.
\ee
The family (\ref{17}) possesses the optimal order four if $p(0)=1,\ p'(0)=2$ and $|p''(0)|<+\infty.$

\smallskip

Proceeding in the same way as in the case of the family of three-point methods (\ref{10}) and using (\ref{9}),
we obtain from (\ref{17}) the fourth order two-point family of iterative methods for finding multiple zeros
    \be
\left\{\begin{array}{l}
y_k=x_k-m \dfrac{f(x_k)}{f'(x_k)},\\[9pt]
x_{k+1}=y_k-m u_kP(u_k) \dfrac{f(x_k)}{f'(x_k)},\quad u_k=\Bigl(\dfrac{f(y_k)}{f(x_k)}\Bigr)^{1/m}, \label{18}
\end{array}
\right.
\ee
where $P(u)$ is the weight function which satisfies
\be
P(0)=1,\ \ P'(0)=2,\label{19}
\ee
that is, its truncated Taylor series is
$$
P(u)=1+2u.
$$
This result is expected since the iterative formula (\ref{18}) coincides with the first two steps of the family (\ref{12}). Explicit expression of $\hat\e=\hat x-\a$ (regarding (\ref{18})) is given by (\ref{14a}) (for $\e_z$).

\smallskip

Several examples of functions that satisfy (\ref{19}) are listed below.
\ba
&&P(u)=\frac{1+\beta u}{1+(\beta-2)u}\ \ \mbox{\rm (of King's type, see \cite{king})},\quad P(u)=\Bigl(1+\frac{2u}{r}\Bigr)^r,\\
&&P(u)=\frac{1+\gamma u^2}{1-2u},\quad P(u)=\frac{1}{1-2u+a u^2},\quad P(u)=\frac{u^2+(c-2)u-1}{cu-1}.
\ea
where $r\in \mbox{\bm{$Q$}}$ and $\gamma,\ a,\ c\in \mbox{\bm{$R$}}$ are arbitrary parameters.

\bigskip

4.6 {\it Equivalence of methods}

    \bigskip

    Some words about the equivalence of the methods presented in Section 2 and 3, and the methods derived by Schr\"oder's approach in Section 4.
    \medskip

    {\bf Equivalence 1.} In Section 2 it was shown that the method (\ref{2}) presented in \cite{lee-AMC} (assuming slight simplification by taking $\la=0$ in (\ref{2})) is a special case of the method (\ref{1}) proposed by Zhou et al. \cite{zhou-AMC}. Comparing  iterative methods (\ref{1}) and (\ref{18}) it is evident that  both formulas are equivalent, which is easily obtained taking $uP(u)\equiv G(u).$ Moreover, the family (\ref{18}) requires only two conditions $P(0)=1,\ P'(0)=2$ compared to three conditions $G(0)=0,\ G'(0)=1,\ G''(0)=4.$ It is important to note that the derivation of the family (\ref{1}) and latter convergence analysis are more complicated than for the family (\ref{18}). Furthermore, the derivation of (\ref{18}), based on Schr\"oder's approach, is entirely natural and crystal clear.
    \medskip

    {\bf Equivalence 2.} The family (\ref{8}), obtained by natural choice $a_1=1,\;a_2=0$ in the family (\ref{4}) (proposed in \cite{behl}), is equivalent to the family (\ref{12}), which is evident setting $S(u)\equiv mP(u)$ and $R(u,v)\equiv mQ(u,v).$ The family (\ref{4}) (and, consequently, (\ref{8})) was presented without derivation and motivation, while the family (\ref{12}) was derived on an easy and obvious way using Schr\"oder's approach.
\bigskip

{\bf Remark 5.} Apart from the family (\ref{3}), other optimal multipoint methods for finding a simple zero (some of them are presented in the book \cite{elzevir}) can be transformed by
introducing $F(x)=f(x)^{1/m}$ to multipoint methods for approximating a multiple zero keeping optimal order of convergence. This subject is left to readers. However, the author of this paper does not expect new papers in this direction since such methods are pretty expensive. This is the subject of the next section.

\section{Root-ratio methods are not competitive}

In  Section 4 we have demonstration a general scheme for constructing multipoint methods for multiple zeros using  basic iterative formulas for simple zeros and Schr\"oder's $f^{1/m}$-idea. We observe that all previously presented methods deal with real or complex values of the forms
$$
\Bigl(\frac{f(y_k)}{f(x_k)}\Bigr)^{1/m},\quad \Bigl(\frac{f(z_k)}{f(x_k)}\Bigr)^{1/m},\quad\Bigl(\frac{f(z_k)}{f(y_k)}\Bigr)^{1/m}.
$$

Computer algebra systems and computation software of digital computers often meet the problem of finding the $m$-th root for arbitrary $m.$ In the case of  specific values of $m$  they find the {\it principal value} of the $m$-root $z^{1/m}$ ($m$ is natural number given as numerical entry) among $m$ values of the sought root in the form
\be
z^{1/m}=|z|^{1/m}\Bigl(\cos\tfrac{\theta}{m}+i\,\sin\tfrac{\theta}{m}\Bigr),\ \  \theta={\rm Arg}\,z\in (-\pi,\pi).\label{20}
\ee

 From (\ref{20}) it is clear that the computation of the $m$-root consumes a lot of CPU time. The following test, implemented on PC with Intel-i7 processor and  clock speed 2.8 GHz, has given  CPU times (expressed in $\mu sec$) for different values of $m$ in calculation of $\bigl((a+ib)/(c+id)\bigr)^{1/m}.$ For the authenticity of the test, one million experiments have been performed taking random numbers for $a,b,c,d$ in each cycle to eliminate the use of possibly memorized data from previous cycles. In the real case we set $b=0,\ d=0.$  The average CPU times for one evaluation are given in Table 1.

 \begin{table}[htb]
 \begin{center}
 {\footnotesize
    \begin{tabular}{|l|l|l|l|l|l|l|l|}
 \hline
   $m$  & 1 & 2 & 3 & 4 & 5 & 6 & 7 \\ \hline
 CPU (in $\mu sec$), real case & 5.25 & 22.7 & 32.24 & 31.86 & 32.25 & 32.03 & 33.1\\
 \hline
 CPU (in $\mu sec$), complex case & 13.77 & 54.82 & 64.42 & 64.07 & 65.83 & 68.05 & 66.3 \\ \hline
   \end{tabular}
 \caption{CPU times in the calculation of the $m$-th root of real and complex numbers\label{t1}}
}
 \end{center}
 \end{table}

 In our experiments  we observed that the CPU times for $m\ge 3$  almost do not change in the real as well as the complex case (see Table 1). In this way we are able to find reliable ratio of computation times for $m\ge 3$ (multiple zeros) and $m=1$ (simple zeros):
  \be
\text{\rm Real case:}\quad \frac{{\rm CPU}_{(m\ge 3)}}{{\rm CPU}_{(m=1)}} \approx 6.2, \qquad
\text{\rm Complex case:}\quad  \frac{{\rm CPU}_{(m\ge 3)}}{{\rm CPU}_{(m=1)}} \approx 4.8. \label{21}
 \ee

\medskip
{\bf Remark 6.} For comparison purpose, we performed one million calculations  of the value of a polynomial of degree 20 with real coefficients and complex argument, both chosen randomly in each cycles (Horner's scheme was used). Average CPU times for one evaluation was 116 $\mu sec$, which is only two times slower than CPU time in calculation of $m$-th root in the case of complex numbers (see Table 1). Therefore, root-ratio methods are expensive form a computational point of view.

\medskip

From (\ref{21}) and Remark 6 we can draw trustworthy conclusion that root-ratio multipoint methods, such as (\ref{1}), (\ref{2}), (\ref{3}), (\ref{4}), (\ref{5}), (\ref{7}), (\ref{8}), (\ref{15b}), (\ref{7b})  and other non-listed methods (if there exist), are pretty inefficient. As mentioned in Remark 6, further work on the construction of root-ratio methods is pointless and does not make any advance in the topic. Combining various weight functions in order to derive ``new methods" is rather a kind of play and inevitably leads to  minor modifications without a proper importance. This fact was emphasized in \cite{AMC-SURVEY} but the construction of modest modifications of  original contributions has continued. Different iterative formula does not mean automatically that an advance in the topic was achieved.

\smallskip

Is there a good alternative method which is efficient and convenient for applications? The answer is yes, and it is very likely known to
many authors who work in the area of iterative processes. Li et al. proposed in \cite{li-2009} the following (optimal) two-point method of order four:
\be
\left\{\begin{array}{l}
y_k=x_k-\frac{2m}{m+2}\cdot\dfrac{f(x_k)}{f'(x_k)},\\[9pt]
x_{k+1}=x_k-\dfrac{\frac{1}{2}m(m-2)A_m t_k-\frac{m^2}{2}}{1-R_m t_k}\cdot\dfrac{f(x_k)}{f'(x_k)},\quad t_k=\dfrac{f'(y_k)}{f'(x_k)} ,\ A_m=\Bigl(\frac{m+2}{m}\Bigr)^m. \label{22}
\end{array}
\right.
\ee
 Zhou et al. \cite{zhou-2011} later proposed the generalization of the method (\ref{22}) in the form
\be
\left\{\begin{array}{l}
y_k=x_k-\frac{2m}{m+2}\cdot\dfrac{f(x_k)}{f'(x_k)},\\[9pt]
x_{k+1}=x_k-\phi(t_k)\cdot\dfrac{f(x_k)}{f'(x_k)},\quad t_k=\dfrac{f'(y_k)}{f'(x_k)}, \label{23}
\end{array}
\right.
\ee
which has the order four under the specific conditions for the weight function $\phi.$  See, also, \cite{li-2010}. The notion of ``convenient" is explained in Remark 8.

\medskip
We conclude this paper with two remarks of general interest.

\bigskip

{\bf Remark 7.} The presence of arbitrary parameters in any zero-finding iterative formula makes sense only if these parameters improve characteristics of  presented methods (such as acceleration of convergence, wider domain of convergence, more accurate approximations, etc.). Otherwise, from an algorithmic point of view, free parameters should be chosen so that an iterative formula is as simple as possible -- numerical analysts and programmers will always choose the simpler formula in such a way that the best characteristics of the employed methods are maintained. Inserting numerous useless parameters does not make a method better or more general in the genuine sense. Unfortunately, many authors construct ``novel" iterative formulas by adding parameters in an artificial way or by varying different weight functions. In essence, such methods are only modest modifications of existing methods and  offer a little  contribution to the topic.

\bigskip

{\bf Remark 8.}
It should be emphasized that very high accuracy of solutions of nonlinear equations, provided by root-solvers of order eight or more, is not needed for solving a huge number of practical problems; fourth order methods (such as (\ref{22}) and (\ref{23}))  produce quite satisfactory results in practice. The question {\it ``how many decimals of zero approximations do we really need in practice?"} is equivalent to the question {\it ``how many decimals of $\pi$  do we really need in practice"?"} A pretty convincing answer can be found  in the issue of {\it NASA/JPL Edu}, March 16, 2016:

\smallskip
$\bullet$ For interplanetary navigation with spacecraft Voyager 1 (launched in 1977, distant from Earth about 22 billion km),  Jet Propulsion Labaratory (California Institute of California, Pasadena, USA) and NASA use very accurate calculations involving $\pi$ with most 15 decimal digits! Not more! The distance error is about 5 cm!

\smallskip
$\bullet$ The radius of the visible universe is about 46 billion light years. To express the circumference of a circle with this radius via the diameter a hydrogen atom (the simplest atom) we need at most 40 decimal digits of $\pi!$

\bigskip
{\bf Acknowledgement.} This work was supported by the Serbian Ministry of Education and Science under Grant 174022.

\thebibliography{10}
\medskip

{

\bibitem{zhou-AMC} X. Zhou, X. Chen, Y. Song, Families of third and fourth order methods for multiple roots of nonlinear equations, Appl. Math. Comput, 219 (2013), 6030--6038.

     \bibitem{lee-AMC} M. Y. Lee, Y. I. Kim, \'A. A. Magren\'an, On the dynamics of tri-parametric family of optimal fourth-order multiple-zero finders with a weight function of the principal mth root of a function-to-function ratio, Appl. Math. Comput. 315 (2017), 564--590.

          \bibitem{liu}  B. Liu, X. Zhou,
A new family of fourth-order methods for multiple roots
of nonlinear equations, Nonlinear Anal. Model. Control, 2013 (18), 143--152.

\bibitem{behl} R. Behl, A. Cordero, S. Motsa, J. R.  Torregrosa,  An eighth-order family of optimal multiple root finders and its dynamics. Numer. Algor. 77 (2018), 1249--1272.

   \bibitem{zafar-1} F. Zafar, A. Cordero, R. Quratulain, J. R. Torregrosa,  Optimal iterative methods for finding multiple roots of nonlinear equations using free
parameters. J. Math. Chem. 56 (2018), 1884--1901.

 \bibitem{zafar-3} F. Zafar, A. Cordero, S. Sultana, J. R. Torregrosa,
Optimal iterative methods for finding multiple roots of nonlinear equations using weight functions and dynamics, J. Comput. Appl. Math. 345 (2018), 352--374.

\bibitem{KT} H. T. Kung, J. F. Traub, Optimal order of one-point and
multipoint iteration, Journal of the ACM  21 (1974), 643--651.

\bibitem{sreder}  {E. Schr\"oder}, {\"Uber unendlich viele Algorithmen zur
Aufl\"osung der Gleichungen},   Math. Ann. 2 (1870), 317--365.

\bibitem{stjuart} {G.W. Stewart}, {On infinitely many algorithms for solving equations}, drum.lib.umd.edu/\linebreak handle/1903/577 (CS-TR-2990.ps) (University of Maryland Library).

 \bibitem{AMC-2011} J. D\v zuni\'c, M. S. Petkovi\'c, L. D. Petkovi\'c, A family of optimal three-point methods for solving nonlinear equations using two parametric functions, Appl. Math. Comput. 217 (2011), 7612--7619.

\bibitem{elzevir} M. S. Petkovi\'c, B. Neta, L. D. Petkovi\'c,
J. D\v zuni\'c, Multipoint Methods for Solving Nonlinear
Equations, Elsevier,
Amsterdam-Boston-Heidelberg-London-New York, 2013.

\bibitem{chun} C. Chun, Some fourth-order iterative methods for solving nonlinear equations, Appl. Math. Comput. 195 (2008), 454--459.

\bibitem{petko-AADM} M. S. Petkovi\'c, L. D. Petkovi\'c, Families of optimal multipoint methods for solving nonlinear equations: a survey, Appl. Anal. Discrete Math. 4 (2010), 1--22.

\bibitem{king} R. F. King, A family of fourth order methods for nonlinear equations, SIAM J. Numer. Anal. 10 (1973), 876--879.

    \bibitem{AMC-SURVEY}   M. S. Petkovi\'c, B. Neta, L. D. Petkovi\'c, J.
D\v zuni\'c,  Multipoint methods for solving nonlinear
equations: a survey, Appl. Math. Comput.
226 (2014) 635--660.

\bibitem{li-2009} S. Li, X. Liao, L. Cheng, A new-fourth-order iterative method for finding multiple root of nonlinear equations, Appl. Math. Comput. 215 (2009), 1288--1292.

    \bibitem{zhou-2011} X. Zhou, X. Chen, Y. Song, Construction of higher order methods for multiple roots of nonlinear equations, J. Comput. Appl. Math. 235 (2011), 4199--4206.

\bibitem{li-2010} S. Li, L. Cheng, B. Neta, Some fourth-order nonlinear solvers with closed formulae for multiple roots, Comp. Math. Appl. 59 (2010), 126--135 .

\bibitem{AMC-2011} J. D\v zuni\'c, M. S. Petkovi\'c, L. D. Petkovi\'c, A family of optimal three-point methods for solving nonlinear equations using two parametric functions, Appl. Math. Comput. 217 (2011) 7612--7619.

}

\end{document}